\documentclass[10pt,twoside]{article}
\usepackage{Latex-document}
\almvol{00}
\almttone{Frontiers of Science Awards}
\almtttwo{for Math (2024)}
\firstpage{1}
\usepackage[english]{babel}
\usepackage[utf8]{inputenc}
\usepackage{amsmath,nicefrac}
\usepackage{mathtools} 
\usepackage{empheq}
\usepackage{multicol,float}

\usepackage{tabu}
\usepackage[nocompress]{cite}

\usepackage{amsmath}
\usepackage{amssymb}

\numberwithin{equation}{section}

\DeclareMathOperator{\Tr}{Tr}
\newcommand{\rmd}{\mathrm{d}}
\newcommand{\R}{\mathbb{R}}
\newcommand{\N}{\mathbb{N}}
\newcommand{\transpose}{^{*}}
\renewcommand{\top}{*}

\let\oldparagraph\paragraph
\renewcommand{\paragraph}[1]{\oldparagraph{#1.}}

\begin{document}

\markboth{\hfill{\rm Jiequn Han, Arnulf Jentzen, \& Weinan E} \hfill}{\hfill {\rm Deep BSDE method \hfill}}

\title{A brief review of the Deep BSDE\\ method for solving high-dimensional partial differential equations}

\author{Jiequn Han, Arnulf Jentzen, Weinan E}

\begin{abstract}
High-dimensional partial differential equations (PDEs) pose significant challenges for numerical computation due to the curse of dimensionality, which limits the applicability of traditional mesh-based methods. Since 2017, the Deep BSDE method~\cite{EHanJentzen2017,HanJentzenE2018} has introduced deep learning techniques that enable the effective solution of nonlinear PDEs in very high dimensions. This innovation has sparked considerable interest in using neural networks for high-dimensional PDEs, making it an active area of research. In this short review, we briefly sketch the Deep BSDE method, its subsequent developments, and future directions for the field.
\end{abstract}

\maketitle

\tableofcontents

\section{Introduction}
Partial differential equations (PDEs) are fundamental tools for modeling a wide range of natural and engineering phenomena. However, developing numerical algorithms to solve high-dimensional PDEs has long been a formidable challenge. Traditional mesh-based methods, which rely on specific choices of basis, inevitably face the notorious ``curse of dimensionality." Consequently, practical algorithms for tackling high-dimensional problems have been very limited. The work~\cite{HanE2016deepcontrol} was the first to adapt the deep learning framework to solve high-dimensional scientific computing problems, with a focus on stochastic control problems, an area closely linked to Hamilton--Jacobi--Bellman (HJB) PDEs. Shortly after, the work~\cite{EHanJentzen2017,HanJentzenE2018} introduced the Deep BSDE method, a practical algorithm capable of solving nonlinear PDEs in very high dimensions (hundreds or even thousands) for a broad class of problems. This method was the first numerical approach based on modern deep learning to effectively address general nonlinear PDEs in high dimensions, opening new possibilities across various disciplines that require solutions to high-dimensional PDEs. Since its introduction, the use of neural networks to solve PDEs has become one of the most active research areas in applied and computational mathematics.

The following sections provide a brief summary of the Deep BSDE method, review the subsequent advances it has inspired, and discuss future prospects for research in high-dimensional PDEs.

\section{Summary of the method}
We consider the following general semilinear parabolic PDE 
\begin{multline}
\label{eq:PDE}
  \partial_t u 
  + \mu(t, x) \cdot \nabla_x u 
  + \tfrac{1}{2} \Tr\bigl(\sigma(t,x) [ \sigma(t,x) ]^\top ( \operatorname{Hess}_x u ) \bigr) 
\\
  + f( t, x, u, [\sigma(t,x)]^{\top} (\nabla_x u) \big) 
  = 0
\end{multline}
for $ t \in (0,T) $, $ x \in \mathbb{R}^d $. Given the terminal condition 
$
  u(T,x ) = g(x),
$
we aim to compute its solution 
$ u \colon [0,T] \times \R^d \to \R $
at time $ t = 0 $. 
Here \( \mu \colon [0, T] \times \mathbb{R}^d \to \mathbb{R}^d \)
and 
\(\sigma \colon [0, T] \times \mathbb{R}^d \to \mathbb{R}^{d \times d} \) are sufficiently regular functions  
and for every $ d \times d $-matrix $ A \in \R^{ d \times d } $, we use
$ A^{ \top } $ 
and $ \Tr( A ) $ to
denote the transpose and the trace of $ A $, respectively. 
The Deep BSDE method approaches this problem by reformulating the PDE as a variational problem, 
particularly a stochastic optimal control problem. This is accomplished through the use of backward stochastic differential euqations (BSDEs) and deep neural networks, which is why it is termed as 
``Deep BSDE method”. 
As a by-product, the method also provides an efficient algorithm for solving high-dimensional BSDEs.

Recall that the classical Feynman--Kac formula (cf.\ \cite{kac1949distributions,oksendal2013stochastic}) 
gives a stochastic representation for solutions of linear PDEs of the Kolmogorov type. BSDE formulations
can be interpreted as a nonlinear generalization to the Feynman--Kac formula 
for a class of nonlinear parabolic PDEs. 
Specifically, we consider the forward stochastic differential equation (SDE)
\begin{equation}
\label{eq:forward_x}
  X_t = \xi + \int_0^t\mu( s, X_s )\,\rmd s +\int_0^t \sigma( s, X_s ) \, \rmd W_s
\end{equation}
for $ t \in [0,T] $
where $\xi$ is a random variable with support in $ D \subseteq \mathbb{R}^d $ 
denoting the region in which we are interested in solving the solution $ u(0,\cdot) $ at time $t=0$. 
When we evaluate the PDE solution along this process, 
using It\^{o}'s lemma and the PDE in \eqref{eq:PDE}, 
we obtain that
\begin{equation}
\label{eq:backward_u}
\begin{split}
  u(t, X_t) - u(r, X_r)
& 
  = - \int_r^t f\big( 
  s, X_s, u(s, X_s), [ \sigma( s, X_s ) ]^{ \top } ( \nabla_x u )( s, X_s )
  \big) \, \rmd s \\
  & \quad + \int_r^t [ \nabla u( s, X_s ) ]^{ \top } \,\sigma( s, X_s )\, \rmd W_s
\end{split}
\end{equation}
for $ r, t \in [0,T] $ with $ r < t $
where we only know the terminal condition $ u(T,X_T) = g(X_T) $.

The structure of the equations in \eqref{eq:forward_x}--\eqref{eq:backward_u} can be summarized abstractly as BSDEs~\cite{Pardoux1992}:
\begin{equation}
\begin{split}
\label{eq:BSDE}
    X_t & = \xi + \int_0^t\mu( s, X_s )\, ds +\int_0^t \sigma( s, X_s ) \, \rmd W_s , 
\\
    Y_t & = g( X_T ) + \int_t^T f( s, X_s, Y_s, Z_s ) \, \rmd s 
    - \int_t^T ( Z_s )^{ \top } \, \rmd W_s.
\end{split}
\end{equation}
It was shown in \cite{Pardoux1992,Pardoux1999} that there is an 
appropriate up-to-equivalence unique stochastic process
$ ( X_t, Y_t, Z_t ) $, $ t \in [0,T] $, 
with values in
$ \mathbb{R}^d \times \mathbb{R} \times \mathbb{R}^d $ 
that satisfies the equations in \eqref{eq:BSDE} above.

Furthermore, the derivation above illustrates the connection between the BSDE in \eqref{eq:BSDE} 
and the PDE in \eqref{eq:PDE}.
Let $ u \colon [0,T] \times \mathbb{R}^d \to\mathbb{R} $ be a solution of the PDE in \eqref{eq:PDE} and define
\begin{equation}
\label{eq:nonlinear_Feynman_Kac}
  Y_t = u( t, X_t )
\qquad  
  \text{and}
\qquad 
  Z_t = [ \sigma( t, X_t ) ]^{ \top }  
  ( \nabla_x u )( t, X_t ).
\end{equation}
for $ t \in [0,T] $. 
Then $ ( X_t, Y_t, Z_t ) $, $ t \in [0,T] $, is a solution for the BSDE in \eqref{eq:BSDE}. 
Conversely, by the uniqueness of the BSDE, solving the BSDE in~\eqref{eq:BSDE} also gives the solution to the PDE in~\eqref{eq:PDE}.

Given the equivalence between the BSDE and the PDE, 
we can also reformulate the original PDE problem as the following variational problem:
\begin{align}
&\inf_{Y_0,\{Z_t\}_{0\le t \le T}} \mathbb{E}\big[ |g(X_T) - Y_T|^2 \big], \\
&s.t.\quad X_t = \xi + \int_{0}^{t}\mu(s,X_s)\, \,\rmd s + \int_{0}^{t}\sigma(s,X_s)\, \rmd W_s, \\
&\hphantom{s.t.}\quad Y_t = Y_0 - \int_{0}^{t}f(s,X_s,Y_s,Z_s)\,  \rmd s + \int_{0}^{t}(Z_s)\transpose\, \rmd W_s.
\end{align}
Note that given $Y_0$ and $\{Z_t\}_{0\le t \le T}$, the process of $Y_t$ above is also defined forwardly. Therefore, one can view this problem as a stochastic control problem, given the process of $X_t$, one chooses starting point $Y_0$ and controls the process $Y_t$ through $Z_t$ such that it matches $g(X_T)$ at the terminal time. The minimizer of this variational problem is the solution to the PDE and vice versa.

\begin{figure}[H]
\centering
\includegraphics[width=0.9\textwidth]{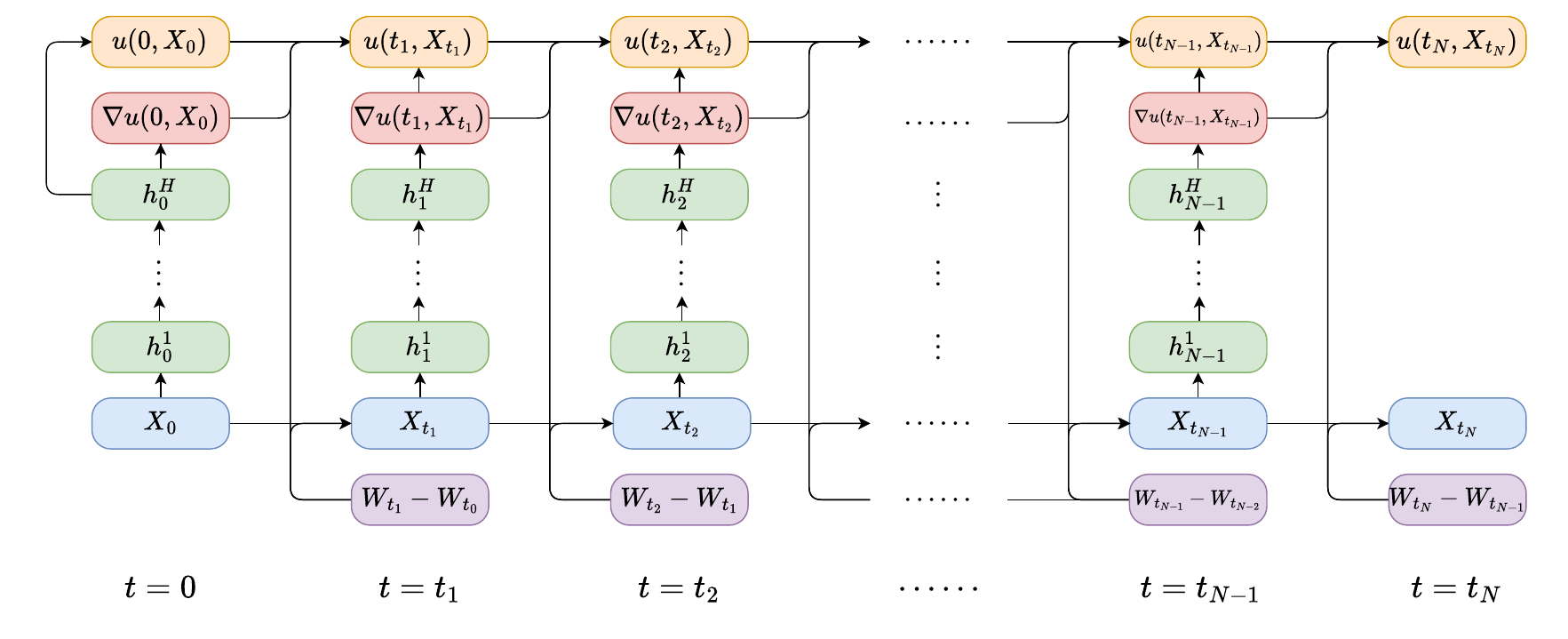}
\caption{Network architecture for solving parabolic PDEs. 
Each column corresponds to a subnetwork at time $ t = t_n $. Reprinted from \cite{EHanJentzen2017}.
}
\label{fig}
\end{figure}

The remaining task is to solve the above variational problem numerically. A key idea of the Deep BSDE method 
is to approximate the unknown functions
$X_0 \mapsto u(0, X_0)$  for $Y_0$ and
$X_t \mapsto [ \sigma(t,X_t) ]^{ \top } ( \nabla_x u )(t,X_t)$ for $Z_t$ (see \eqref{eq:nonlinear_Feynman_Kac})
by feedforward neural networks $\psi$ and $\phi$. Then, we discretize time, say, by using the Euler 
scheme on a grid $0 = t_0 < t_1 < \ldots < t_N = T$:
    \begin{align}
    &\inf_{\psi_0, \{\phi_n\}_{n=0}^{N-1}} \mathbb{E} |g(X_T) - Y_T|^2, \\
    &s.t.\quad X_0 = \xi, \quad Y_0 = \psi_0(\xi), \\
    &\hphantom{s.t.}\quad X_{t_{n+1}} = X_{t_n} + \mu(t_n,X_{t_n})\Delta t_n + \sigma(t_n,X_{t_n})\Delta W_n, \\
    &\hphantom{s.t.}\quad Z_{t_n} = \phi_n(X_{t_n}), \\
    &\hphantom{s.t.}\quad Y_{t_{n+1}} = Y_{t_n} - f(t_n,X_{t_n},Y_{t_n},Z_{t_n})\Delta t_n + (Z_{t_n})\transpose\Delta W_n. \label{eq:discretized_Y}
    \end{align}
 At each time step $t_n$, we associate a subnetwork.
By stacking these subnetworks together, we form a deep neural network. This network takes the paths 
$\{ X_{ t_n } \}_{ 0 \leq n \leq N }$ and 
$\{ W_{ t_n } \}_{ 0 \leq n \leq N }$ 
as the input data and produces the final output, denoted by
    $\hat{u}( 
    \{ { X_{ t_n } } \}_{ 0 \leq n \leq N } , 
    \{ W_{ t_n } \}_{ 0 \leq n \leq N } )$, 
as an approximation to 
$u( t_N, X_{ t_N } )$. 
The network architecture is illustrated in Figure~\ref{fig}, 
which is reprinted from the original work \cite{EHanJentzen2017}.
The resulting network has a very natural ``residual neural network'' structure~\cite{he2016deep} embedded 
in the stochastic difference equations in \eqref{eq:discretized_Y}.

The error in matching the given terminal condition defines the loss function:
  \begin{equation}
  \label{eq:BSDE_objective}
    l(\theta) = 
    \mathbb{E}\Big[
      \big|g( X_{ t_N } ) - \hat{u}\big(\{ X_{ t_n } \}_{ 0 \leq n \leq N } , \{ W_{ t_n } \}_{ 0 \leq n \leq N }\big)\big|^2
    \Big].
  \end{equation}
From a machine learning perspective, this loss function does not require the generation of training data in advance. The initial condition $X_{t_0}$ and the paths $\{ W_{ t_n } \}_{ 0 \leq n \leq N }$ serve as the data, and they are generated on the fly very easily. As a result, this model can be thought of as having access to an infinite amount of data, making it particularly suitable for training with stochastic gradient descent (SGD). Note that we denote the different networks by $\psi_0$ and $\{\phi_n\}_{n=0}^{N-1}$. These networks can have independent or shared parameters, which does not impact the application of SGD. Additionally, suppose we are only interested in the PDE solution at a single spatial point $x_0$. In that case, we can set $\xi = \xi_0$ as a deterministic value, simplifying $\phi_0 = u(0, \xi_0)$ and $\psi_0 = \nabla_x u(0, \xi_0)$ to deterministic values that can also be optimized.

By leveraging stochastic reformulation and deep neural network approximation, 
the Deep BSDE method proposed above can solve several high-dimensional PDEs in the region of interest efficiently. We refer to~\cite{EHanJentzen2017,HanJentzenE2018} for various numerical examples.

\section{Subsequent advances}
\label{sec:subsequent_advances}
The success of the Deep BSDE method has sparked a wide range of research efforts aimed at leveraging deep learning to tackle high-dimensional PDEs, 
leading to significant advancements in various scientific and engineering fields. In this section, we briefly review selected subsequent works that 
focus on solving high-dimensional PDEs, particularly those pose intrinsic challenges for classical methods, 
and we also briefly mention selected research efforts regarding the theoretical analysis 
of such deep learning-based approximation methods for high-dimensional PDEs. For more comprehensive discussions and further references, 
we refer, for example, to the review papers~\cite{e2021algorithms,germain2021neural,Beck_2023,Blechschmidt}.

\paragraph{Methods based on BSDEs and related stochastic Feyman--Kac formulations} The efficiency 
and applicability of BSDE-based methods have been enhanced in various ways; see, for example, \cite{beck2019machine,fujii2019asymptotic,nusken2021solving,hure2020deep,ji2020three,pham2021neural,nusken2023interpolating,han2024deeppicard}. In particular, the stochastic optimization objective in~\eqref{eq:BSDE_objective} derived from BSDEs 
is well-suited to the prevailing deep learning optimization paradigm based on SGD. This framework can be naturally extended to various contexts, 
including the actor-critic algorithm for the HJB equation~\cite{zhou2021actor,ata2023drift}, 
Nash equilibrium in multi-agent or mean-field games~\cite{Han2019deep,han2024learning,lauriere2022learning}, 
elliptic problems~\cite{Han2020solving,han2020derivative,nusken2023interpolating}, 
as well as seconds order BSDEs and fully-nonlinear PDEs~\cite{beck2019machine,han2024deeppicard}. 
We also refer, for example, to \cite{beck2021deep,MR4293960,BernerDablanderGrohs2020arXiv}
for deep learning approximation methods for (parametric) linear and semilinear Kolmogorov PDEs 
based on classical linear Feynman--Kac formulations using standard forward SDEs.

\paragraph{Methods based on the least-squares formulation}
Besides stochastic reformulations based on BSDEs and SDEs, 
there are other widely used formulations for designing numerical algorithms for PDEs, among which the least-squares formulation is one of the most versatile. This approach transforms a PDE problem into a variational problem by minimizing the squared residual of the PDE. In classical numerical analysis, this method is often less preferred because the resulting numerical problems tend to be more poorly conditioned compared to those from other formulations. By replacing the space of trial functions with the hypothesis space of neural networks, 
this formulation has led to the development of physics-informed neural networks (PINNs) \cite{Sirignano2018dgm,raissi2019physics} 
(although the method in \cite{Sirignano2018dgm} is typically referred to as Deep Galerkin method). Notably, such deep learning methods for PDEs based on the least-squares formulation were proposed and studied as early as the 1990s in the context of low-dimensional PDEs (see, for example, \cite{lagaris1998artificial,dissanayake1994neural,jianyu2003numerical}).
Originally developed for solving low-dimensional PDEs, PINNs have quickly gained popularity due to their generality and ease of implementation. 
Recently, they have also been adapted to address high-dimensional PDEs when combined with Hutchinson trace estimation~\cite{hu2024hutchinson}. 
For further references in this line of research we also refer, for example, to the survey articles \cite{CuomoEtAl2022,karniadakis2021physicsinformedML}.

\paragraph{Methods based on the Ritz formulation}
The Ritz formulation is a classical variational approach for solving PDEs. By using neural networks as trial functions, this approach naturally leads to the development of the Deep Ritz method in \cite{E2018deep}. A closely related area is variational Monte Carlo, which is used to solve for the eigenfunctions of the Schr\"{o}dinger equation. There has been growing interest in constructing neural network-based wavefunctions for this problem, starting with the influential work~\cite{carleo2017solving} on spin systems, and later extending to many-electron systems~\cite{han2019solving,luo2019backflow, pfau2019ab, Hermann2019deep}.

\paragraph{Methods based on the Galerkin formulation}
The Galerkin formulation is another popular approach for solving PDEs in engineering and mathematical modeling, based on the weak formulation of a PDE that involves both trial and test functions. Designing machine learning-based algorithms using the Galerkin formulation is much less straightforward than the previous formulations, as it is not a standard optimization-based approach. The work~\cite{Zang2020weak} innovatively interprets the weak solution as the solution of a min-max problem and then develops an algorithm to train neural network-based approximations to PDEs in an adversarial way (cf., for example, also \cite{chen2023friedrichs,ValsecchiOliva2022xnodewan}).

\paragraph{Mathematical analysis of deep learning-based approximation methods for PDEs}
The rapid numerical progress in this field has opened new avenues for theoretical investigations for deep neural network-based approximations for high-dimensional PDEs. While some advances have been made (see, for example, \cite{Han2018convergence,hutzenthaler2020proof,Knochenhaueretal2022arXiv,lu2021priori,chen2021representation,han2022convergence,
duan2022convergence,grohs2023proof,carmona2021convergence,carmona2022convergence,MR4127967}), a full understanding of the overall error in deep learning methods—including approximation, generalization, and optimization errors—remains incomplete, even in the relatively simpler case of one-dimensional settings. The approximation error that arises when deep neural networks are used to approximate solutions of PDEs, along with the generalization error to some extent, appears to be comparatively better understood from a theoretical perspective. Particularly, several results in the literature rigorously prove that deep neural networks can approximate solutions of high-dimensional PDEs without the curse of dimensionality (cf., for example, \cite{MR4127967,Bezneaetal2024arXiv,cheridito2023efficient,ElbraechterSchwab2018,MR4454916,grohs2023proof,MR4283528} and the references therein for such approximation results for linear PDEs and cf., for instance, \cite{Ackermannetal2024arXiv,Cioica-Lichtetal2022arXiv,hutzenthaler2020proof,neufeld2024rectifiedGradientDep,neufeld2023pide} and the references therein for such approximation results for semilinear PDEs with Lipschitz nonlinearities). Specifically, such results reveal that the number of real parameters of the approximating deep neural networks grows only at most polynomially in both the reciprocal $ \nicefrac{ 1 }{ \varepsilon } $ of the prescribed approximation accuracy $ \varepsilon > 0 $ and the PDE dimension $ d \in \N $.

\section{Future prospects}

As the field of machine learning-based methods for solving high-dimensional PDEs continues to evolve, several areas with significant applications stand out as promising directions for future research. While existing works have made considerable progress, there remain challenges that require further development to fully realize the potential of these approaches. Below, we highlight a few key areas.

\paragraph{Optimal Control and HJB Equation} Optimal control problems often involve solving high-dimensional HJB equations, where the curse of dimensionality was originally coined and is closely related to the BSDE formulation discussed earlier. Recent works~\cite{nakamura2021adaptive,zhang2022initial} have shown promising results using deep learning techniques for optimal feedback control. However, the complexity of real-world problems still demands further advancements. In financial mathematics and economics, optimal decision-making in high-dimensional and stochastic environments is a common challenge, and extensive research has been conducted in these directions~\cite{becker2020pricing,liang2021deep,beck2021deep,gnoatto2023deep,payne2024deep,domingo2024stochastic}.

\paragraph{Probabilistic Modeling and High-Dimensional Densities} Probabilistic modeling involving high-dimensional densities also faces the curse of dimensionality. In recent years, learning methods based on PDEs, such as the Fokker-Planck equation, and SDEs have achieved significant advances in solving and sampling high-dimensional densities~\cite{boffi2024deep,song2021scorebased,han2024deeppicard}, in a similar spirit to the Deep BSDE method introduced earlier. However, many challenges remain in improving the accuracy and efficiency of these methods, making this a rich area for future research.

\paragraph{Quantum Mechanics and the Schr\"{o}dinger Equation} The Schr\"{o}dinger equation is fundamental to quantum mechanics, yet finding accurate solutions in high dimensions remains a significant challenge. While neural network-based approaches, as mentioned earlier, have shown potential in approximating wavefunctions with high precision, issues such as scalability and the need for more efficient training algorithms are still unresolved. Additionally, the time-dependent Schr\"{o}dinger equation has been much less explored compared to its time-independent counterpart, presenting another avenue for future research.

\paragraph{Plasma Physics and Kinetic Equations} Kinetic equations, such as the Boltzmann equation, are crucial for modeling the behavior of gases and plasmas. Although the dimensionality of these equations is not exceedingly high, the complexity of particle physics results in intimidating computational costs. Since the work~\cite{han2019uniformly}, deep learning-based closures for kinetic equations have emerged as a promising approach to reduce simulation costs~\cite{huang2022machine,li2023learning}. Future work could include improving the scalability of these methods to larger-scale problems and enhancing the interpretability of the resulting models.

\paragraph{Theoretical Advancements}
Despite the progress that has been made in the theoretical analysis of deep learning methods for high-dimensional PDEs, developing a complete convergence analysis remains a fundamental open problem for essentially all deep learning-based approximation methods for PDEs (cf.\ \cite{gonon2023random}). While several open questions remain regarding approximation and generalization errors—particularly for high-dimensional PDEs with gradient-dependent nonlinearities, such as HJB equations (cf.\ \cite{HutzenthalerJentzenKruse201912,neufeld2024rectifiedGradientDep})—the optimization error of deep learning methods poses an especially significant challenge. This difficulty persists even in the case of one-dimensional PDEs and one-dimensional target functions (cf., for instance, \cite{gentile2022approximation,ibragimov2022convergence,welper2023approximation}).

\subsection*{Acknowledgements}
The second author gratefully acknowledges partial support from the National Science Foundation of China (NSFC) under grant number 12250610192. In addition, the second author gratefully acknowledges the Cluster of Excellence EXC 2044-390685587, Mathematics Münster: Dynamics-Geometry-Structure funded by the Deutsche Forschungsgemeinschaft (DFG, German Research Foundation). The third author gratefully acknowledges support from the National Natural Science Foundation of China (NSFC) under grant number 92270001.

\bibliographystyle{acm}
\bibliography{ref}

\address{(JH) Center for Computational Mathematics, Flatiron Institute, New York, USA.\\
\email{jhan@flatironinstitute.org}}

\address{(AJ) School of Data Science and Shenzhen Research Institute of Big Data, The Chinese University of Hong Kong, Shenzhen (CUHK-Shenzhen), China; \\
Applied Mathematics: Institute for Analysis and Numerics, University of M\"{u}nster, Germany.\\
\email{ajentzen@cuhk.edu.cn, ajentzen@uni-muenster.de}}

\address{(WE) Center for Machine Learning Research and School of Mathematical Sciences, Peking University, Beijing, China;\\
AI for Science Institute, Beijing, China.\\
\email{weinan@math.pku.edu.cn}}
\end{document}